\documentclass[11pt]{amsart}
\usepackage{latexsym}
\usepackage{amsfonts}
\usepackage{amsmath}
\usepackage{amssymb}
\usepackage{lcy}
\usepackage{verbatim}
\newtheorem{Th}{Theorem}[section] \newtheorem{Cor}[Th]{Corollary}
\newtheorem{Lem}[Th]{Lemma} 
\newtheorem{Claim}[Th]{Claim}
\newtheorem{Rem}[Th]{Remark}

\newtheorem{claim-num}{Claim}

\numberwithin{equation}{section}
\renewcommand{\theequation}{\thesection.\arabic{equation}}

\def\aut#1{\operatorname{Aut}(#1)}
\def\inn#1{\operatorname{Inn}(#1)}

\def\inv{^{-1}}

\def\str#1{\langle#1\rangle}

\def\e{\varepsilon}
\def\eps{\varepsilon}

\def\N{\mathbf N}
\def\Z{\mathbf Z}

\renewcommand{\le}{\leqslant}
\renewcommand{\ge}{\geqslant}

\def\wid{\operatorname{wid}}

\def\wg{\operatorname{WG}}
\def\bu{\mbox{\boldmath{$u$}}}
\def\SL#1{\operatorname{sl}(#1)}
\def\sgn{\operatorname{sign}}

\begin{document}

\title{On the palindromic and primitive widths of a free group}
\author{Valery Bardakov }
\address{Institute of Mathematics\\
Siberian Branch Russian Academy of Science\\
630090 Novosibirsk\\
Russia}
\thanks{V. Bardakov was partially supported
by the RFFI grant \# 02-01-01118}
\email{bardakov@math.nsc.ru}
\author{Vladimir Shpilrain}
\address{Department of Mathematics\\ The City College of New York\\ New York, NY 10031}
\email{shpil@groups.sci.ccny.cuny.edu}
\author{Vladimir Tolstykh}
\address{Department of Mathematics\\ Kemerovo State University\\ Kemerovo\\  Russia}
\thanks{V. Tolstykh was supported
by a NATO PC-B grant via The
Scientific and Technical Research Council of Turkey (T\"UBITAK)}
\email{tvlaa@mail.ru}
\maketitle

\section*{Introduction}

Let $G$ be a group and $S \subseteq G$ a subset
 that generates $G.$ For each $x \in G$ define
the {\it length} $l_S(x)$ of $x$ relative to $S$ to be
the minimal $k$ such that $x$ is a product of $k$
elements of $S.$ The supremum of the values $l_S(x)$, $x \in G$,
is called the {\it width} of $G$ with respect to $S$
and is denoted by $\wid(G,S).$ In particular,
$\wid(G,S)$ is either a natural number or $\infty.$
If $\wid(G,S)$ is a natural number, then
every element of $G$ is a product of at most $\wid(G,S)$ elements of $S.$

Many group-theoretic results can be interpreted as
determining the width of a group with respect to one 
generating set or another. 
We can mention results on the width of matrix groups
with respect to  the set of transvections 
(see e.g. \cite{Bar:Mat,CarKel:Mat,VasWhel1:Mat})
and the study of the width of verbal subgroups of various free
constructions (\cite{Bar:VS,Dobr:VS,Faiz:VS,Smir}).

It is also worth mentioning that  the
concept of the width of a group can be useful in 
 model theory.  Possible
applications here are based on the following simple argument.
Consider a first-order definable subset $D$ of a
group $G$; in other words, $D$ is the set of
realizations of a first-order formula
in the language of groups. Then the subgroup $\str D$ generated by $D$ is
also definable in $G,$ provided that the width of
$\str D$ relative to $D$ is finite. For example,
one of the authors proved that the family
of all inner automorphisms determined
by powers of primitive elements is definable
in the automorphism group $\aut{F_n}$ of the 
 free group $F_n$ of rank $n \ge 2$ \cite{ToJLM}. It is
however an open question whether or not 
the subgroup $\inn{F_n}$ of all inner
automorphisms is  definable in $\aut{F_n}.$
An affirmative answer would follow from the finiteness
of the width of $F_n$ relative to the set of
all primitive elements (the {\it primitive width}
of $F_n.$) Similarly, inner automorphisms from $\aut{F_n}$
determined by palindromic words in $F_n$ with respect
to a fixed basis of $F_n$ are also definable -- this
time with suitable definable parameters -- in
$\aut{F_n}.$ Palindromic words generate $F_n,$ and 
this again raises the question of  whether or not  the
corresponding ({\it palindromic}) width of $F_n$ is
finite.

The  goal of the present paper is to determine the
primitive and  palindromic widths of a free group.
The main results of the paper show that both widths
of a finitely generated non-abelian free group are
infinite, and for infinitely generated free groups,
the palindromic width is   infinite, too, whereas
the primitive width is not (one easily checks that
the primitive width of an infinitely generated free
group is two).  Note that, although our results do not
solve the problem of definability of the subgroup of
inner automorphisms in the automorphism group of a
finitely generated free non-abelian group, they can
be considered an evidence of  its difficulty.

In Section 1 we examine the palindromic
width of a free group. Recall that given a
basis $X$ of a free group $F$, a reduced word $w \in
F$ is called a {\it palindrome} if $w$ reads the
same  left-to-right and  right-to-left as a word in 
letters from $X^{\pm 1}.$
The {\it palindromic width} of $F$ is therefore the
width of $F$ with respect to the palindromes
associated with a given basis of $F.$ Palindromes of free
groups have already proved useful in  studying 
various aspects of combinatorial group
theory: for instance, in
\cite{Coll,GloJen}, palindromic automorphisms of free
groups are studied and in  \cite{Helling}, palindromes
are used in a   description of automorphisms of a 
two-generator free group. 

Using  standard methods developed earlier for the
study of verbal subgroups of free contructions, we
prove that the palindromic width of any free group is
infinite (Theorem \ref{Pal-Width-of-a-Free-Group}).
Then we  discuss  relations between
the primitive and palindromic widths of a
two-generator free group $F_2.$ In this special case
we establish the following result which seems to be of
independent interest: any primitive element of $F_2$
is a product of at most two palindromes.
(Here we mention, in passing,  an interesting related fact \cite{ShYu}:
any primitive element of a free associative algebra of rank 2 is palindromic.)
This, together with
Theorem \ref{Pal-Width-of-a-Free-Group} implies
 that the primitive width of $F_2$ is infinite.

In Section 2 we generalize the latter result by proving
that the primitive width of a finitely generated
free group $F_n, n \ge 2,$ is infinite 
(Theorem \ref{Prim-Width-of-a-Free-Group}).
Actually, a stronger result is established: the width of $F_n$ relative to
the set $W_n$ of elements of $F_n$ whose Whitehead graph has
a cut vertex is infinite. Recall that, according
to  a classical result of Whitehead, the set
$W_n$ contains all primitive elements of $F_n$
\cite{White}.

 To conclude the Introduction, we mention several interesting,
in our opinion, problems motivated by the results of this paper.
\medskip

\noindent {\bf Problem 1.} Let $F_n$ be a free
group of finite rank $n \ge 2.$ Is there an algorithm to
  determine the primitive length of an element $w$ of $F_n$? 
\medskip

\noindent {\bf Problem 2.} The same question for the palindromic length of an
element of $F_n.$
\medskip

 There are two interesting related problems that can  be formulated in a 
different language: 
 is the Dehn function relative to the set of primitive (resp. palindromic)
elements recursive? Here, by somewhat abusing the language, by
the Dehn function relative to a generating set $M$ of a group $G$ with
a fixed generating set $S$ we mean
a function $f_{S,  M}(n)$ equal to the maximum   length $l_M(w)$ over all
$w \in G$ such that $l_S(w)=n$. With respect to Problems 1 and 2,
 we naturally assume that $S$ is an arbitrary but fixed basis of $F_n$.
 If these Dehn functions turn out to be recursive, it would be
interesting to find out what they are. It is clear that positive answers 
to Problems 1 and 2 would imply that the corresponding Dehn functions 
are recursive. 

 We mention here  that Grigorchuk and  Kurchanov  \cite{GrigKur}
have reported an algorithm to determine the length of an element $w$ of $F_n$
with respect to the set of all conjugates of elements of a {\it fixed basis}
of $F_n$.

 Finally, we note that the primitive width can be defined not only for 
free groups, but for all relatively free groups as well.  In
particular, it is easy to see that  the
primitive width of any free abelian group and any free
nilponent group is two. Smirnova \cite{Smir} proved
that the primitive width of any free metabelian group
of rank $\ge 2$ is at most four.  Lapshina \cite{Lap}
generalized this to  relatively free groups
with nilpotent commutator subgroup. Later she also proved that
the primitive width of the free solvable group of rank $n$
and class $s$ is at most $4+(s-1)(n-1)$. The same is true,
in fact, for any free polynilpotent group of rank $n$
and class $(k_1, ..., k_s)$ \cite{Lap1}. All this motivates the following
\medskip

\noindent {\bf Problem 3.} Describe relatively free groups whose  primitive
width is infinite.

\medskip

The authors are grateful to  Oleg Belegradek and Oleg
Bogopolsky for helpful discussions.

\section{The palindromic width of a free group}

Let $F$ be a free group with a basis $X.$ We call a
word $w \in F$ a {\it palindrome in letters $X$} if
$$
w = x_{i_1} x{i_2} \ldots x{i_k} = x{i_k} \ldots x{i_2} x_{i_1},
$$
where $x_{i_j} \in X^{\pm 1}$ and $x_{i_j} x_{i_{j+1}} \ne 1$ for
all $j=1,2,\ldots,k-1.$ Equivalently, it can be
said that $w \in F$ is a
palindrome if and only if $w$ coincides with its {\it
reverse} word \cite{Coll}.  For example, if
$x_1,x_2,x_3$ are distinct elements of $X$, then
the word
$$
x_1^{-2} x_2^{3} x_3^{-4}x_2^{3} x_1^{-2}
$$
is a palindrome.

Since  elements of a basis $X$ of a free group $F$
are palindromes in letters $X,$ the set of all palindromes
generates $F$, and we can define the {\it palindromic width}  of $F$ 
as the width relative to the set of palindromes.

\begin{Th} \label{Pal-Width-of-a-Free-Group}
The palindromic width of a non-abelian
free group $F$ is infinite.
\end{Th}

\begin{proof} A standard   idea of showing 
that the width of a given group relative to a given 
  generating set is infinite is 
  constructing a {\it quasi-ho\-mo\-mor\-phism}
$\Delta : F \to \Z$ satisfying the inequality
\begin{equation}
\Delta(uw) \le \Delta(u)+\Delta(w)+\text{const}
\end{equation}
for all $u,w \in F.$ After constructing
such a quasi-homomorphism, we prove, for
all $k \in \N$, that there is a bound $c_k$
such that
$$
\Delta(w) \le c_k
$$
for all $w \in F$ that are products of at most $k$
palindromes. This will surely hold if the values of
$\Delta$ with (\theequation) are reasonably bounded on 
palindromes (say, $\Delta(p)=0$ for all palindromes
$p$); one can say then, somewhat informally,  that $\Delta$
 {\it recognizes} palindromes.  Finally, we shall 
find a sequence $(w_n)$ of elements of $F$ with
$$
\lim_{n \to \infty} \Delta(w_n)=+\infty.
$$
Clearly, having all these requirements on $\Delta$
satisfied, one obtains that the palindromic width of
$F$ is infinite, since for all $n$ with $\Delta(w_n) >
c_k,$ the word $w_n$ is not a product of $k$
palindromes.

 We shall now define  a
required quasi-homomorphism by induction on the number
of syllables of a reduced word $w \in F$
in letters $X.$

1) Suppose  $w$ has exactly one syllable, then
$$
\Delta(w) =0.
$$

2) Let now
$$
w = v_1\ldots v_{m-1}v_m,
$$
where $m> 1$ and $v_k$ are syllables of $w.$ Assume that 
$$
v_{m-1} =a^{\pm k}, v_m =b^{\pm l}, 
$$
where $a,b \in X$, $a \ne b$, 
and $k,l \in \N.$ Then
$$
\Delta(w) =\Delta(v_1\ldots v_{m-1})+ \sgn(l-k).
$$
The function $\sgn : \Z \to \{-1,0,1\}$
takes the value $1$ at positive integers,
the value $-1$ at negative integers and $\sgn(0)=0.$

For example, if $x_1,x_2,x_3$ are distinct
elements of $X,$ then
\begin{align*}
\Delta(x_1^2 x_2^{-3} x_3^4 x_2^{-3} x_1^2) &=
\sgn(2-3) + \sgn(3-4)+
\sgn(4-3)+\sgn(3-2)\\
&=-1-1+1+1=0.
\end{align*}

One readily verifies the following properties
of $\Delta:$

\begin{Claim} \label{Bas-Props-of-f}
{\em (i)} For any word $w$ in $F$, 
$$
\Delta(w)+\Delta(w\inv)=0;
$$

{\em (ii)} for every palindrome $p$ in letters
$X$, 
$$
\Delta(p)=0.
$$
\end{Claim}

Note that not only palindromes give the value of
$\Delta$ equal to $0.$ For instance,
$$
\Delta(x_1 x_2 \ldots x_n ) =0, 
$$
where $x_i$ are distinct elements of $X$.

Now we are going to show  that $\Delta$ is a quasi-homomorphism.

\begin{Lem} \label{Delta1-Is-a-Quasi-Hom}
For any $u,w \in F$
$$
\Delta(u w) \le \Delta(u)+\Delta(w)+6.
$$
\end{Lem}

\begin{proof}
Suppose that
\begin{align*}
u &= u_0 t^k z\inv,\\
w &= z t^m w_0,
\end{align*}
where
\begin{itemize}
\item $u_0,w_0$ are reduced;
\item $t \in X$ and $k+m \ne 0$;
\item neither $u_0$ ends with $t,$ nor
$w_0$ begins with $t.$
\end{itemize}
Then we have
\begin{align*}
\Delta(u) &=\Delta(u_0) +\eps_1+\eps_2+\Delta(z^{-1}),\\
\Delta(w) &=\Delta(z)+\eps_3+\eps_4+\Delta(w_0),
\end{align*}
where $\eps_1,\ldots,\eps_4 \in \{-1,0,1\}.$
Furthermore, 
$$
\Delta(uw) =\Delta(u_0 t^{k+m} w_0) =\Delta(u_0)+\eps_5 +\eps_6+\Delta(w_0), 
$$
and again $\eps_5,\eps_6 \in \{-1,0,1\}.$
Clearly,
$$
\eps_5+\eps_6 \le \eps_1+\eps_2+\eps_3+\eps_4+6.
$$
Hence by Claim \ref{Bas-Props-of-f}
$$
\Delta(u w) \le \Delta(u)+\Delta(w)+6,
$$
as  needed.
\end{proof}

If now $w$ is a product of at most $k$ palindromes,
then a straightforward induction argument using
Lemma \ref{Delta1-Is-a-Quasi-Hom} proves that
$$
\Delta(w) \le 6k-6.
$$
We are going to find a sequence $(w_n)$ 
of words of $F$ such that
$$
\Delta(w_n)=n-1
$$
for all $n \ge 1.$  This will imply that for any
$n > 6k-5$ the word $w_n$ is not
a product of $k$ palindromes,  and the palindromic
width of $F$ is therefore infinite. A particular choice of $w_n$ 
like that is as  follows :
$$
w_n =x_1 x_2 x_1^2  x_2^2\ldots x_1^n x_2^n,
$$
where $n \ge 1$ and $x_1,x_2$ are fixed distinct elements of $X.$

This completes the proof of  Theorem \ref{Pal-Width-of-a-Free-Group}.
\end{proof}

In the remainder of this section we discuss 
the relation between the palindromic
and primitive widths of two-generator
free groups. We start with an elementary observation
from \cite{GloJen} which provides a useful
description of palindromes.

\begin{Lem} \label{Pals-As-Theta-Invs} Let $F$ be a free group with a basis $X.$
Suppose that $\theta$ is an involution of $\aut F$
that inverts each element $x \in X.$ Then $w \in F$ is
a palindrome in letters $X$ if and only if $\theta$
inverts $w.$
\end{Lem}

\begin{Rem} \label{My-Invs}
\em Note that if $p$ is a palindrome in letters $X$ and
$u$ is any word of $F(X),$ then
by Lemma \ref{Pals-As-Theta-Invs} the words $\theta(u)u\inv$
and $\theta(u) p u\inv$ are also
palindromes in letters $X$
(cf. \cite[Lemma 2.3]{ToJLM}.)
\end{Rem}

\begin{Lem} \label{Prim-As-Two-Pals}
Any primitive element of a two-generator free group
$F_2$ is a product of at most two palindromes.
\end{Lem}

\begin{proof}
Let  $\{x,y\}$ be a basis of $F_2.$
Let $\tau_w$ denote the inner automorphism
determined by $w \in F_2.$ Suppose that $\theta$ is the automorphism of $F_2$ 
that inverts both $x$ and $y.$ Since $\theta$ induces the
automorphism $\text{id}_A$ on the abelianization $A=F_2/[F_2,F_2]$
of $F_2$,  any product of the form
$$
\theta \sigma \theta \sigma\inv,
$$
where $\sigma$ is an arbitrary automorphism
of $F_2$, induces the trivial automorphism
of $A$ and hence, by  a well known result of Nielsen
(see e.g. \cite[I.4.13]{LSch}), it is an inner
automorphism determined by some $p \in F_2$:
$$
\theta \sigma \theta \sigma\inv =\tau_p.
$$
We claim that $p$ is a palindrome.  Indeed, 
\begin{align*}
\theta \tau_p \theta &= \tau_{\theta(p)}
                     =\theta (\theta \sigma \theta \sigma\inv) \theta \\
             &=\sigma \theta \sigma\inv\theta
             =\tau_p\inv=\tau_{p\inv}.
\end{align*}
Then
$\theta(p)=p\inv$ and by Lemma
\ref{Pals-As-Theta-Invs} $p$ is a palindrome; we
denote $p$ by $p(\sigma).$

Let $U$ be the Nielsen automorphism that
acts on $\{x,y\}$ as follows:
\begin{align*}
U(x) &=x,\\
U(y) &=x y.
\end{align*}

One easily checks that
$$
\theta U \theta U\inv=\tau_{x}.
$$

Then we have
\begin{align*}
\tau_{p(\sigma U)} &=\theta \sigma U \theta U\inv \sigma\inv =
\theta \sigma \theta (\theta U \theta U\inv) \sigma\inv \\
&=\theta \sigma \theta \tau_{x} \sigma\inv =
\theta \sigma \theta \sigma\inv \sigma \tau_{x} \sigma\inv \\
&=\tau_{p(\sigma)} \tau_{\sigma(x)}.
\end{align*}
This tells us that
$$
\sigma(x) = p(\sigma)\inv p(\sigma U),
$$
that is, any element in the orbit
of $x$ under $\aut{F_2}$ (in other words,
any primitive element of $F_2$) is a product of two
palindromes in letters $x,y.$
\end{proof}

As a corollary we immediately obtain
the following fact which will be
generalized in the next section. By the {\it primitive width} of a free
group $F$ we mean the width of $F$ relative 
to the set of all primitive elements.

\begin{Cor} \label{Pal-Width-Then-Prim-Width}
The primitive width of a two-generator
free group $F_2$ is infinite.
\end{Cor}

\begin{proof}
The result follows immediately from 
 Theorem \ref{Pal-Width-of-a-Free-Group}
and Lemma \ref{Prim-As-Two-Pals}.
\end{proof}

\begin{Rem} \label{Some-Rems-on-Pals}
\em (i) It is also possible to derive
a fact similar to Lemma \ref{Prim-As-Two-Pals}
from one of the results in the paper \cite{Helling} by Helling.
Indeed, let $a$ be a primitive element of $F_2=F(x,y).$
It then follows from the Theorem on p. 613 of \cite{Helling} that there are a
palindrome $p$ in letters $x,y$ and an element $z \in F_2$ such that
$$
a =z y\inv p x z\inv
$$
or
$$
a =z x\inv p y z\inv.
$$
On the other hand, for instance, in the
first case we have
$$
a=z y\inv \theta(z\inv) \cdot
\theta(z) p z\inv \cdot
z x \theta(z\inv) \cdot
\theta(z) z\inv, 
$$
where $\theta$ is the automorphism of $F_2$ that
inverts both $x$ and $y.$ Therefore, by Remark
\ref{My-Invs}, an arbitrary primitive element of $F_2$
is  a product of at most four palindromes.

(ii) There is an important consequence
of the cited result by Helling which
seems to pass unnoticed  and which
improves the description of primitive
elements of $F_2$ given in  
\cite{CoMeZi}. It is proved there 
that the conjugacy class of every primitive
element of $F_2(x,y)$ contains either
an element of the form
\begin{equation}
a=x y^{m_1} x y^{m_2} \ldots x y^{m_k},
\end{equation}
where $m_k \in \{n, n+1 \}$ for a suitable fixed
natural $n,$ or an element obtained
from $a$ by (a) permutation of the basis
elements $x,y$, or (b) by inversion
of one or both  basis elements. It is then clear in view of the result by Helling
that the distribution of exponents
$n,n+1$ in (\theequation) must
obey a rather strict rule: there
must be a cyclic permutation $a'$
of the word $a$ such that $x a' y\inv$
or $x\inv a' y$ is a palindrome.

(iii) In  sharp contrast with the two-generator case, 
 the palindromic length of primitive elements
in a free group $F_n$ of rank   $n > 2$ {\it
cannot} be uniformly bounded. Indeed, if there is a
word $w_k(x,y)$ in letters $x,y$ which is not a
product of at most $k$ palindromes in $F_2=F(x,y),$
then, for instance, the word $z w_k(x,y),$ a primitive
element of $F_3=F(x,y,z),$ cannot be written as a
product of at most $k$ palindromes. It follows immediately
from the fact that the homomorphism
$$
x \to x, y \to y, z \to 1
$$
from $F(x,y,z)$ to $F(x,y)$ takes
palindromes to palindromes.
\end{Rem}

 To  conclude this section, we discuss a possible
generalization of our results here. 

The notion of a palindromic word, or
a palindrome, can be defined
for free products of groups as follows. Given a free
product
\begin{equation} \tag{$*$}
G = \left.\prod_{i \in I}\right.^* G_i,
\end{equation}
a reduced word $g$ of $G$ is called
a {\it palindrome} associated with
the decomposition $(*)$ if, when written in syllables from $G_i,$
$g$ ``reads the same left-to-right and  right-to-left".
The palindromic width of $G$ is the
width of $G$ relative to the set of  palindromes.

It seems that any free product of groups which is not a 
 product of two cyclic groups of order two has
infinite palindromic width.  We  show this here 
in the case of a free product of two {\it infinite} groups. 
Using the ``homomorphism argument", as in Remark
\ref{Some-Rems-on-Pals} (iii), one can see
 that any free product of groups at least two of which
are infinite has infinite palindromic width.
The proof is, in fact, a modification of the proof of
Theorem \ref{Pal-Width-of-a-Free-Group}.

Suppose that  $G= A*B$ and both groups $A,B$ are
infinite. Then we consider a surjective map $d : A \cup B\to \N$ 
such that $d(A)=d(B)=\N.$ We define a
quasi-homomorphism $\Delta$ recognizing palindromes as
follows:

1) $\Delta(v)=0,$ if $v \in A,B;$

2) if $g=v_1\ldots v_n$ is a reduced
word written in syllables and $n > 1,$ then
$$
\Delta(g) =\Delta(v_1\ldots v_{n-1}) +\sgn( d(v_n)-d(v_{n-1})).
$$
One then verifies that $\Delta$ is indeed a quasi-homomorphism.
Furthermore, let $(a_n)$ and $(b_n)$ be
sequences of elements of $A$ and $B$ respectively
with $d(a_n)=d(b_n)=n$ for all $n \in \N.$
Then
$$
\Delta( a_1 b_1 a_2 b_2 \ldots a_n b_n) =n-1.
$$
The reader will easily supply missing details
in our sketch.

\section{The primitive width of a free group}

\begin{Th} \label{Prim-Width-of-a-Free-Group}
Let $F$ be a free group. Then the primitive
width of $F$ is infinite if and only if $F$
is finitely generated. The primitive
width of any infinitely generated free group is two.
\end{Th}

\begin{Rem}
\em 
In contrast, the palindromic width of an infinitely
generated free group is infinite, as we have seen above. Thus
the set of palindromes in any infinitely generated free
group appeares to be  quite ``sparse" compared to the 
set of primitive elements.  \end{Rem}

\begin{proof}
We begin with the proof of our second, 
almost obvious, statement. Suppose that $X$ is an infinite basis of $F.$ Take an
arbitrary $w$ in $F$ and assume that $w$ is a word in
letters $y_1,\ldots,y_k$ of $X$:
$w=w(y_1,\ldots,y_k).$ Since $X$ is infinite,
there is a basis letter $x$ which is different
from all $y_1,\ldots,y_k.$ We have then
$$
w = x\inv \cdot x w(y_1,\ldots,y_k).
$$
Clearly, both words $x\inv$ and $x w(y_1,\ldots,y_k)$
are primitive elements of $F.$

Let now $F=F_n$ be a finitely generated
free group with a free basis $X=\{x_1,\ldots,x_n\}.$
 Let $w$ be a (reduced) word
from $F.$ The {\it Whitehead graph} $\wg_w$ of $w$
is constructed as follows: the
vertices of $\wg_w$ are the elements
of the set $X^{\pm 1}$, and the edges
are determined by the pairs $(a,b\inv),$
where $a,b \in X^{\pm 1}$ are such that
there is an occurence of the subword
$ab$ in $w.$ Note that, like in  most  recent
papers that use the Whitehead graph (see e.g. \cite{BuVe}),
 we are in fact using a  simplified 
version of the graph introduced by Whitehead in
\cite{White}.
The principal result about the 
Whitehead graph of a primitive element of a 
free group is the following:

\begin{Th}[Whitehead, \cite{White}]
The Whitehead graph of a primitive
element of a free group has a cut
vertex.
\end{Th}

Recall that a vertex $v$ of a graph $\Gamma$ is said to
be a {\it cut vertex} if  removing 
the vertex $v$ along with all its adjacent edges from $\Gamma$ 
increases the number of connected components of the graph. 
 The fact that the above theorem
is also true for the simplified Whitehead graph 
is an easy corollary of Theorem 2.4 from the paper
\cite{St} by Stallings.

Now we give a proof of the existence of words of $F$ that
  cannot occur as proper subwords of primitive
elements of $F.$  We claim that any word $u \in F$
whose Whitehead graph is Hamiltonian, that is, such
that all the vertices of $\wg_u$ can be included into a
simple circuit, meets our condition. Indeed, if $u$
occurs in a word $w$ as a proper subword, then the
graph $\wg_u$ is a subgraph of $\wg_w,$ and hence
$\wg_w$ contains no cut vertices since there are no
cut vertices in $\wg_u$ (because the latter graph is Hamiltonian). 
Thus, by the Whitehead theorem, 
$w$ is not primitive.

An example of such a word $u$ can be easily
found: for instance, the
Whitehead graph of the word
$$
\bu = (x_1^2 \ldots x_n^2)^2
$$
is Hamiltonian.

\begin{Lem}
Let $k$ be a natural number, and $\bu = (x_1^2 \ldots x_n^2)^2$.
 The element $\bu^{2k}$ is not a product
of at most $k$ primitive elements of $F.$
\end{Lem}

\begin{proof}
Let $\SL{w}$ denote the syllabic
length of a word $w \in F$ with
respect to the basis $X$ of $F.$
Suppose that
$$
\bu^{2k} = p_1 \ldots p_m,
$$
where $m \le k$ and $p_1,\ldots,p_m$
are primitive elements of $F.$ After
possible reductions the product $p_1\ldots p_m$
turns into the product $p_1'\ldots p_n',$
where $p_i'$ is a subword of $p_i,$
and the latter product has the following
property:
\stepcounter{equation}
\begin{quote}
\item[(\theequation)] for every $i=1,\ldots,m-1$ the last
syllable of $p_i'$ and the first syllable of $p_{i+1}'$
do not annihilate each other; in other words,
$$
\SL{p_i'p_{i+1}'}=\SL{p'_i}+\SL{p'_{i+1}}+\e_i,
$$
where $\e_i$ equals $0$ or $-1.$
\end{quote}

We then have that
$$
\SL{\bu^{2k}} = \SL{p_1'\ldots p'_m},
$$
or
$$
4nk = \SL{p_1'}+\eps_1+\ldots +\eps_{m-1}+\SL{p'_m},
$$
where $\eps_i=0,-1.$ Letting
$l$ denote the sum $-\eps_1-\ldots-\eps_{m-1},$
we get
$$
4nk+l = \SL{p_1'}+\ldots +\SL{p'_m}.
$$
Then there is an index  $j$ such that
$$
\SL{p'_j} \ge \frac{4nk+l}{m} \ge \frac{4nk+l}{k} \ge 4n \ge 2n+2.
$$
The condition (\theequation) implies therefore
that $p'_j$ contains at least $2n$ consecutive
syllables of $\bu.$ However, for the subword $z$
of $\bu$ determined by these syllables
the Whitehead graph is clearly Hamiltonian, because
$z$ is  the result of
a cyclic permutation of the syllables
of $\bu.$ It follows
that $p_j$ has $z$ as a proper subword,
and then it is not primitive.  This contradiction completes the proof
of the lemma.
\end{proof}

According to Lemma \theLem, there are
elements of $F$ of arbitrarily large
primitive length. Therefore the primitive
length of $F$ is infinite.
\end{proof}

\end{document}